\documentclass[a4paper,12pt]{article}
\usepackage[english]{babel}
\usepackage{amssymb}
\usepackage{amsmath}
\usepackage{amsthm}
\usepackage{graphicx, epsfig, color}
\newtheorem{lemma}{Lemma}[section]
\newtheorem{theorem}{Theorem}[section]
\author{Dmitry Efimov\thanks{email: defimov@dm.komisc.ru}}
\date{Department of Mathematics,\\ 
Komi Science Centre UrD RAS,\\
Syktyvkar, Russia}
\title{Determinants of generalized binary band matrices\thanks{The study is supported by Program of UD RAS, project 15-16-1-3}}

\begin{document}
\maketitle
\abstract{Under binary matrices we mean matrices whose entries take one of two values.
In this paper, explicit formulae for calculating  the determinant of some type of binary Toeplitz matrices are obtained.
Examples of the application of the determinant of binary Toeplitz matrices for the enumeration of even and odd permutations
of different types are given.}




\section*{Introduction}
Under binary matrices we will mean matrices whose elements can take only two values.
Matrices of this kind arise in different mathematical questions.
For example, this type includes such popular objects in mathematics and its applications as matrices over the field $GF(2)$ \cite{Boston}. 
The Hadamard's problem of finding the maximal 
determinant of $(-1,1)$-matrices, i.e. matrices consisting from $1$ and $-1$, is well known \cite{Hadamard}.  
Further, $(0,1)$-matrices are one of the favorite objects of the enumerative combinatorics [3--6]. 
Also one can use binary $(1,x)$-matrices for enumerative problems \cite{Shev1}.
These examples can be extended.

One of the basic notion of the matrix theory is the notion of determinant.
There exist effective algorithms for the determinant calculation, 
for example, the modified method of Gaussian elimination, 
which run in polynomial time.
Nevertheless, in cases of some kinds of matrices it is possible to obtain
good explicit formulae for the determinant expression.
On the one hand, these formulae allow to draw certain conclusions about properties of matrices.
On the other hand, they give even greater gain in the speed of the determinant calculation.
Such formulae are known only for very limited class of matrices.

Our paper is motivated to some extent by the recent work \cite{Krav} in which
the explicit formula for the determinant of some binary circulant matrices has been obtained.
In the present paper we get explicit formulae for determinants of some kinds of binary Toeplitz matrices.
Considered matrices are close  in their structure to band matrices,
therefore we call them generalized band matrices.
For our purposes we use quite elementary methods.
Applying the Laplace expansion we obtain recurrent formulae leading to required result.

The paper is organized as follows.
In the second section we give and prove the formulae, 
which allow  to efficiently calculate determinants of 
generalized binary band matrices in an explicit form. 
In the third section we give a few examples of the application of the determinant of such 
binary matrices for the enumeration of even and odd permutations of different types.

\section{The main part}
Let $n, k, l$ be integers, $1 \leq l\leq k\leq n$.
Consider a binary Toeplitz matrix $A=(a_{ij})$ of order $n$ with the following elements:
$$
a_{ij}=\left\{
 \begin{array}{ll}
   b,&\textrm{if } -l<j-i<k;\\
   a,&\textrm{otherwise,}
 \end{array}
\right.
$$
where $a$ and $b$ ($a\not=b$) belong to a commutative associative ring with a unit.
In the case $a=0$ we get a so-called \textit{band matrix} (\cite{Golub}, p.16).
So we can say that we consider \textit{generalized binary band matrices}.
Our purpose is to get explicit formulae for calculating the determinant of the matrix $A$.
We will divide this problem into two cases.
\paragraph{Case $1$. } Let $l=1$. 
In other words, let the first row of the matrix $A$ have the form:
$$
[\underbrace{\overbrace{b\ \dots\ b}^k\ a\  \dots\ a}_n],
$$
and the $(i+1)$-st row is obtained from the first row  by the removal  of $i$ elements on the right
and addition  of $i$ elements $a$ on the left:
\begin{equation}\label{15.06_2}
A=\left(
\begin{array}{ccccccc}
b&\dots&b&&&&\\
&b&\dots&b&&a&\\
&&b&\dots&b&&\\
&&&\ddots&\ddots&\ddots&\\
&&&&b&\dots&b\\
&a&&&&\ddots&\vdots\\
&&&&&&b
\end{array}
\right).
\end{equation}
\begin{theorem}\label{16.06_1}
Let $n\equiv p\pmod{k}$, where $0<p\leq k$.
Then the determinant of the matrix (\ref{15.06_2}) is equal to:
\begin{equation}\label{09.06_1}
  \det{A}=(b-a)^{n-1}\left(b+\frac{n-p}{k}a\right).
\end{equation}
\end{theorem}

For the proof of Theorem \ref{16.06_1} one can apply the method given in \cite{Krav} 
for the proof of explicit formulae of the determinants of circulant matrices 
(using the formula of the determinant of a $2\times 2$ block matrix).
But in this case another method will be more convenient.
First let us formulate some auxiliary statements.
Consider a square matrix of order $n$ of the form:
\begin{equation}\label{15.06_3}
\left(
\begin{array}{ccccccc|c}
b&\dots&b&&&&&a\\
&b&\dots&b&&a&&a\\
&&b&\dots&b&&&a\\
&&&\ddots&\ddots&\ddots&&\vdots\\
&&&&b&\dots&b&a\\
&a&&&&\ddots&\vdots&\vdots\\
&&&&&&b&a\\
\hline
a&a&a&\dots&a&\dots&a&a
\end{array}
\right).
\end{equation}
The given matrix is obtained from the matrix (\ref{15.06_2}) of order $(n-1)$
by addition of one more (the last) row and one more (the last) column, 
consisting entirely of elements $a$.
Let $f_n$ denote the determinant of such matrix.
\begin{lemma}
The determinant of the matrix (\ref{15.06_3}) is equal to:
\begin{equation}\label{15.06_5}
f_n=(b-a)^{n-1}a.
\end{equation}
\end{lemma}

\noindent{\it Proof.}
Subtracting the penultimate row from the last row and expanding the determinant along the last row, 
we get that $f_n=(b-a)f_{n-1}$. Using the given formula consecutively to $f_{n-1}$, $f_{n-2}$ and so on, we get (\ref{15.06_5}). \hfill$\square$\hspace{-1mm}
\vspace{2mm}

Suppose that $k=n$, i.e. consider the matrix of order $n$ of the form:
\begin{equation}\label{15.06_4}
\left(
\begin{array}{ccccc}
b&b&b&\dots&b\\
&b&b&\dots&b\\
&&b&\dots&b\\
&a&&\ddots&\vdots\\
&&&&b
\end{array}
\right).
\end{equation}
By $g_n$ denote the determinant of such matrix.

\begin{lemma}\label{16.06_3}
One can calculate the determinant of the matrix (\ref{15.06_4}) by the formula:
\begin{equation}\label{16.06_2}
g_n=(b-a)^{n-1}b.
\end{equation}
\end{lemma}

\noindent{\it Proof.}
Similarly to the proof of the previous lemma. \hfill $\square$\hspace{-1mm}
\vspace{2mm}

Consider the matrix (\ref{15.06_2}) of order $n$ with $k<n$.
Let $d_n$ denote the determinant of such matrix.
\begin{lemma}
The following equality holds:
\begin{equation}\label{15.06_1}
 d_n=
 \left\{
       \begin{array}{ll}
         (b-a)^kd_{n-k}+(b-a)^{n-1}a,&k<\frac{n}{2};\\
         (b-a)^{n-1}(b+a),&k\geq\frac{n}{2}. 
        \end{array}
     \right.
\end{equation}
\end{lemma}

\noindent{\it Proof.}
Let $k<\frac{n}{2}$, i.e. the number of elements $b$ in the first row of the matrix (\ref{15.06_2}) is less than the number of elements $a$.  
Let us subtract the penultimate row of the matrix from the last row and expand the determinant along the last row.
Performing this procedure consecutively $k$ times, we come to the equality:
$$
d_n=(b-a)^kd_{n-k}+(b-a)^kf_{n-k}.
$$
Substituting here the formula (\ref{15.06_5}), we get the first row of the equality (\ref{15.06_1}).  

Let $\frac{n}{2}\leq k$, i.e. the number of elements $b$ in the first row of the matrix (\ref{15.06_2}) is more or equal to the number of elements $a$, 
but less than $n$.  
Executing the same $k$ consecutive expansions along the last row as in the first case, we come to the equality:
$$
d_n=(b-a)^k(f_{n-k}+g_{n-k}).
$$
Substituting here the formulae (\ref{15.06_5}) and (\ref{16.06_2}), we get the second row of (\ref{15.06_1}).   \hfill $\square$\hspace{-1mm}
\vspace{2mm}

\noindent{\it Proof.}
\textbf{(of the Theorem \ref{16.06_1})}.
Suppose that  $n\equiv p\pmod{k}$, where $0<p\leq k$, i.e. $n=km+p$, where $m$ is a non-negative integer.
Assume that $k<n$.
Then $m>0$ and, using the first row of the formula (\ref{15.06_1}) recurrently $m-1$ times, we obtain the equality:
$$
d_n=(b-a)^{n-k-p}d_{k+p}+(m-1)(b-a)^{n-1}a.
$$  
Since $\frac{k+p}{2}\leq k$, then applying the second row of the formula (\ref{15.06_1}) to $d_{k+p}$
and performing obvious transformations, we get (\ref{09.06_1}).

If $k=n$, i.e. the first row of the matrix (\ref{15.06_2}) entirely consists from elements $b$, 
then $p=n$ and the formula (\ref{09.06_1}) gives us $d_n=(b-a)^{n-1}b$, which corresponds to the statement of the Lemma \ref{16.06_3}. \hfill $\square$\hspace{-1mm}
\vspace{2mm}

\paragraph{Case 2. } Let $l>1$.
In other words, let the $i$-th row of the matrix $A$, $i=1,\dots, l$, have the form:
$$
[\underbrace{\overbrace{b\ \dots\ b}^{k+i-1}\ a\  \dots\ a}_n],
$$
and the $(l+j)$-th row is obtained from the  $l$-th one by  the removal  of $j$ elements on the right 
and addition  of $j$ elements $a$ on the left:
\begin{equation}\label{29.06_1}
A=\left(
\begin{array}{cccccc}
b&\dots&b&&&a\\
\vdots&\ddots&&\ddots&&\\
b&&\ddots&&\ddots&\\
&\ddots&&\ddots&&b\\
&&\ddots&&\ddots&\vdots\\
a&&&b&\dots&b
\end{array}
\right).
\end{equation}
The matrix $A$ is not symmetric in general, but it is persymmetric 
(i.e. symmetric with respect of the secondary diagonal) like all Toeplitz matrices.
\begin{lemma}\label{16.08_1}
In the matrix (\ref{29.06_1}) exactly $\max{(k+l-n,0)}$ rows consist entirely from elements $b$. 
\end{lemma}

\noindent{\it Proof.}
First note that each row of the matrix $A$ contains no more than $k+l-1$ elements $b$ by definition.
Let $k+l-n\leq 0$. Hence $k+l-1<n$, 
i.e.  there is at least one element $a$ in each row and there are no rows, consisting entirely from elements  $b$.
Let now $k+l-n>0$ or, alternatively, $k+l-1\geq n$.
If $n=k+i-1$, $1\leq i\leq l$, then the number of rows,
that consist entirely from elements $b$, will be equal to  $l-i+1=k+l-n$.  \hfill $\square$\hspace{-1mm}
\vspace{2mm}

\begin{theorem}
Let $1<l\leq k\leq n$ and $n\equiv p\pmod{k+l-1}$, $0\leq p < k+l-1$.
Then the determinant of the matrix (\ref{29.06_1}) is equal to:
\begin{equation}\label{09.08_2}
\det{A}=\left\{
\begin{array}{ll}
(-1)^\frac{(k-1)(l-1)n}{k+l-1}(b-a)^{n-1}\left(b+\frac{n-k-l+1}{k+l-1}a\right),&p=0,\\
(-1)^\frac{(k-1)(l-1)(n-1)}{k+l-1} (b-a)^{n-1}\left(b+\frac{n-1}{k+l-1}a\right),&p=1,\\
0,&\text{otherwise}.
\end{array}
\right. 
\end{equation}
\end{theorem}

\noindent{\it Proof.}
By the theorem condition $n\geq 2$.
Let $n\equiv p\pmod{k+l-1}$, $0\leq p < k+l-1$.
This is equivalent to $n=(k+l-1)s+p$, where $s$ is a non-negative integer.

First consider the case when $s=0$, i.e. when $p=n$ and, respectively, $n<k+l-1$.
By lemma \ref{16.08_1}  at least $2$ rows of the matrix $A$ will  consist entirely from elements $b$ in this case, 
therefore the determinant of the matrix $A$ will be equal to $0$, that corresponds to the formula (\ref{09.08_2}). 

Now assume that $s\geq 1$. The first row of the matrix $A$ differs from the second one only by an element in the $(k+1)$-st column.
Let us subtract the first row from the second one and expand the determinant along the second row.
We will get that $\det{A}=(-1)^{k-1}(b-a)\det{A'}$, 
where $A'$ is a matrix, whose first and  second rows differ also only in the $(k+1)$-st column.
Performing this procedure $l-1$ times, we obtain that
$\det{A}=(-1)^{(k-1)(l-1)}(b-a)^{l-1}\det{A''}$, where
$$
A''=\left(
\begin{array}{ccc|cccccc}
b&\dots&b&&&&&&\\
&\ddots&\vdots&\ddots&&&a&&\\
a&&b&\dots&b&&&&\\
\hline
\vphantom{Y^2}&&&b&\dots&b&&&a\\
&&&\vdots&\ddots&&\ddots&&\\
&a&&b&&\ddots&&\ddots&\\
&&&&\ddots&&\ddots&&b\\
&&&&&\ddots&&\ddots&\vdots\\
&&&a&&&b&\dots&b
\end{array}
\right).
$$
The part of the matrix $A''$ that is located above the horizontal line
consists from $k$ rows, and in the first row there are exactly $k$ elements $b$. 
In the lower right corner of the matrix $A''$ there is a submatrix, 
which is formed  by intersection of the last $n-k-l+1$ rows and columns of $A$. 
The given submatrix, in turn,  also has the form (\ref{29.06_1}).

The algorithm can be repeated $s$ times with an obvious shift at each step in the $k$ rows down and in the $k$ columns to the right.
Thus in the second step we subtract the $(k+1)$-st row of the matrix $A''$ from the $(k+2)$-nd one and expand the determinant along the $(k+2)$-nd row,
and repeat this procedure $l-1$ times and so on.
As a result we will get that 
\begin{equation}\label{30.06_3}
\det{A}=(-1)^{(k-1)(l-1)s}(b-a)^{(l-1)s}\det{A'''},
\end{equation}
where
$$
A'''=\left(
\begin{array}{ccccc|c}
b&\dots&b&&a&\\
&\ddots&&\ddots&&a\\
&&\ddots&&b&\\
&a&&\ddots&\vdots&\ddots\\
&&&&b&\dots\\
\hline
&&a&&&M
\end{array}
\right).
$$
The part of the matrix $A'''$ that is located above the horizontal line
consists from $ks$ rows, and in the first row there are exactly $k$ elements $b$. 
In the lower right corner of the matrix $A'''$ there is a submatrix $M$, 
which is formed  by intersection  of the last $n-(k+l-1)s=p$ rows and columns of the matrix $A$. 

It is not hard to see that if  $1<p<k+l-1$ then $M$ contains two identical rows,
hence the matrix $A'''$ also contains two identical rows, and therefore
its determinant and the determinant of the matrix $A$ are equal to $0$. 

If $p=0$ then the matrix $A'''$ is the matrix of the form (\ref{15.06_2}) of order $ks$.
Then by Theorem \ref{16.06_1} we get:
\begin{equation*}
 \det{A'''}=(b-a)^{ks-1}\left(b+\frac{ks-k}{k}a\right)=(b-a)^{ks-1}\left[b+(s-1)a\right].
\end{equation*} 
Then, substituting the given formula in (\ref{30.06_3}) and taking into account that $n=(k+l-1)s$, we finally obtain:
\begin{multline*}
  \det{A}=(-1)^{(k-1)(l-1)s}(b-a)^{(k+l-1)s-1}\left[b+(s-1)a\right]=\\
	=(-1)^\frac{(k-1)(l-1)n}{k+l-1}(b-a)^{n-1}\left(b+\frac{n-k-l+1}{k+l-1}a\right).
\end{multline*}

If $p=1$ then the matrix $A'''$ will be the matrix of the form (\ref{15.06_2}) of order $ks+1$.
Then from Theorem \ref{16.06_1} we get:
\begin{equation*}
 \det{A'''}=(b-a)^{ks}\left(b+\frac{ks+1-1}{k}a\right)=(b-a)^{ks}(b+sa).
\end{equation*}
Then, substituting the given formula in (\ref{30.06_3}) and taking into account that $n=(k+l-1)s+1$, we finally obtain:
\begin{multline*}
 \det{A}=(-1)^{(k-1)(l-1)s}(b-a)^{(k+l-1)s}(b+sa)=\\
=(-1)^\frac{(k-1)(l-1)(n-1)}{k+l-1}(b-a)^{n-1}\left(b+\frac{n-1}{k+l-1}a\right). 
\end{multline*}
\hfill $\square$\hspace{-1mm}
\vspace{2mm}

\section{The application of determinants of binary matrices to the enumeration of permutations}
As mentioned in the Introduction, the binary matrices are one of the favorite objects of the enumerative combinatorics. 
In particular, they are applied for  enumeration of permutations with restricted positions. 
Following \cite{Shev1} let us describe briefly this mechanism. 

Let $A=(a_{ij})$ be a $(0,1)$-matrix of order $n$.
Each of such matrices defines a class $\mathcal{B}(A)$ of restricted permutations.
Namely a permutation $p$ belongs to $\mathcal{B}(A)$ if and only if
the inequality $M_p\leq A$ holds for its incidence matrix $M_p$,
i.e. each element of the matrix $M_p$ is not more than the corresponding element of the matrix $A$.
The matrix $A$ is called the characteristic matrix of the class $\mathcal{B}(A)$.
It is not hard to see that the number of permutations in the class $\mathcal{B}(A)$ is equal to the permanent of the matrix $A$:
$|\mathcal{B}(A)|=\mathrm{per}\, A$.
Denote  the number of even and odd permutations from the class $\mathcal{B}(A)$ by $E_A$ and $O_A$, respectively.
It is obvious that $E_A+O_A=\mathrm{per}\, A$.
It is easy to see also that $E_A-O_A=\det{A}$. 
This implies  the following formulae for calculating the total number
of even and odd permutations from the class $\mathcal{B}(A)$:
\begin{equation}\label{09.08_1}
 E_A=\frac{\mathrm{per}\, A+\det{A}}{2},\ \ \ O_A=\frac{\mathrm{per}\, A-\det{A}}{2}.
\end{equation}

As an example we calculate the number of even and odd permutations $\pi\in S_n$ such that $\pi(i)\not=i,i+1$ for $i=1,\dots n-1$ and $\pi(n)\not=n$. 
The characteristic matrix of this class of permutations is the following $(0,1)$-matrix $A_n$ of order $n$:  
$$
A_n=\left(
\begin{array}{ccccc}
0&0&&&\\
&0&0&1&\\
&&\ddots&\ddots&\\
&1&&0&0\\
&&&&0
\end{array}
\right).
$$
Such matrices arise in the variation of the famous  \textit{m\'enage problem},
where not a round table,  but one side of a rectangular table is considered  (\cite{Raiz}, ch. 8). 

If we denote the permanent of such matrix by $p_n$, then the sequence of permanents 
will satisfy the following recurrence relation:
\begin{equation*}
 (n-1)p_n=(n^2-n-1)p_{n-1}+np_{n-2}+2(-1)^{n+1},\ \ \ p_1=p_2=0. 
\end{equation*}
One can also calculate these permanents  by the following explicit formula:
$$
p_n=\sum_{k=0}^n \binom{2n-k}k (n-k)!(-1)^k. 
$$

Let $d_n$ denote the determinant of the matrix $A_n$. 
Substituting $b=0$, $a=1$, $k=2$ to the formula (\ref{09.06_1}), we get
$$
d_n=(-1)^{n-1}\frac{n-p}{2},\ \ \ n\equiv p\hspace{-2mm}\pmod{2},\ \ \ 0<p\leq 2.
$$
One can also rewrite this equality  in the following form:
$$
d_n=(-1)^{n-1}\left\lfloor\frac{n-1}{2}\right\rfloor=\left\{
\begin{array}{ll}
  \frac{n-1}{2},&\textrm{if } n \textrm{ -- odd};\\
  -\frac{n-2}{2},&\textrm{if } n \textrm{ -- even}.
\end{array}
\right.
$$
Applying formulae (\ref{09.08_1}), we obtain sequences of the number of even ($e_n$) and odd ($o_n$) permutations
of the given type depending on the permutation order:
$$
\begin{array}{c|c|c|c|c|c|c|c|c|c|c|c}
 n&1&2&3&4&5&6&7&8&9&10&\dots\\
 \hline
 p_n&0&0&1&3&16&96&675&5413&48800&488592&\dots\\
 \hline
d_n&0&0&1&-1&2&-2&3&-3&4&-4&\dots\\
\hline
e_n&0&0&1&1&9&47&339&2705&24402&244294&\dots\\
\hline
o_n&0&0&0&2&7&49&336&2708&24398&244298&\dots
\end{array} 
$$

As the second example, we will calculate the number of even and odd permutations $\pi\in S_n$ such that $|\pi(i)-i|>1$, $i=1,\dots,n$. 
The characteristic matrix of the given class of permutations is the matrix of order $n$ $B_n$
whose main diagonal and its neighboring diagonals are zero, and all other elements are equal to $1$:
$$
B_n=\left(
\begin{array}{cccccc}
0&0&&&&\\
0&0&0&&1&\\
&0&0&0&&\\
&&\ddots&\ddots&\ddots&\\
&1&&0&0&0\\
&&&&0&0
\end{array}
\right).
$$
This example is linked with another variation of the m\'enage problem,
where a rectangular table is considered and additional restriction on the placement of men is imposed.
The explicit formula of the total number of such permutations of order $n$ or, alternatively,
of the value of the permanent of 
$B_n$ was found by V.S. Shevelev (see the review \cite{Shev1}). 
Here we will not give it, but indicate only that the sequence $\{p'_n\}$ of such numbers has the id-number $A001883$ in \cite{Sloane}.
Let $d'_n$ denote  the determinant of the matrix $B_n$.
Substituting $b=0$, $a=1$, $k=l=2$ to the formula (\ref{09.08_2}), we get
$$
d'_n=\left\{
\begin{array}{ll}
  \frac{3-n}{3},&\textrm{if } p=0;\\
  \frac{n-1}{3},&\textrm{if } p=1;\\
  0,&\textrm{if } p=2;\\
\end{array}
\right.
$$
where $n\equiv p\pmod{3}$.
Applying formulae (\ref{09.08_1}), we obtain sequences of the number of even ($e'_n$) and odd ($o'_n$) permutations
of the given type depending on the permutation order:
$$
\begin{array}{c|c|c|c|c|c|c|c|c|c|c|c}
 n&1&2&3&4&5&6&7&8&9&10&\dots\\
 \hline
 p'_n&0&0&0&1&4&29&206&1708&15702&159737&\dots\\
 \hline
 d'_n&0&0&0&1&0&-1&2&0&-2&3&\dots\\
 \hline
 e'_n&0&0&0&1&2&14&104&854&7850&79870&\dots\\
 \hline
 o'_n&0&0&0&0&2&15&102&854&7852&79867&\dots
\end{array} 
$$

Another version of the application of binary matrices to enumeration of permutations has been described also in \cite{Shev1}. 
It consists in the following.
Let us consider  a binary matrix $A$ of order $n$ in which some elements are equal to the variable $b$, and the other elements are equal to $1$.  
It is easy to see that the coefficient on $b^k$ in the permanent  $\textrm{per}\, A$ will be equal to the number of permutations $\pi\in S_n$
whose incidence matrices have exactly $k$ units in the positions, in which there are elements $b$ in the matrix $A$.
Respectively, the coefficient on $b^k$ in the determinant $\det{A}$
will be equal to the difference between the number of even and odd permutations of such type.
Then calculating in expressions $\frac{1}{2}(\det{A}\pm \textrm{per}\, A)$ the coefficient on $b^k$, 
we get the number of even and odd permutations of such kind.
 
Let us give an example of the application of this scheme. Let $\pi$ be a permutation. 
Recall that a number $i$ for which $\pi(i)\geq i$ is called  \textit{a weak excedance} of $\pi$  (\cite{Stan}, p. 40).
Let us find the number of even and odd permutations of order $n$ with exactly $k$ weak excedances.
Consider a binary matrix of order $n$ $C_n$ whose elements on the main diagonal and above it are equal to $b$, 
and other elements are equal to $1$:
$$
C_n=\left(
\begin{array}{ccccc}
b&b&b&\dots&b\\
&b&b&\dots&b\\
&&b&\dots&b\\
&1&&\ddots&\vdots\\
&&&&b
\end{array}
\right).
$$
It is well known (\cite{Stan}, p. 39--40) that the permanent of this matrix is equal to the Eulerian polynomial of order $n$,
and, respectively, the coefficient on $b^k$ is equal to the Eulerian number $T(n,k)$ (A008292 in \cite{Sloane}).
By the formula  (\ref{09.06_1}) we get that $\det{C_n}=(b-1)^{n-1}b$.
It follows that the coefficient on $b^k$ in $\det{C_n}$ is equal to:
$$
c(n,k)=(-1)^{n-k}\binom{n-1}{k-1}.
$$
Hence the number of even and odd permutations of order $n$ with exactly $k$ weak excedances is equal to:
\begin{equation*}
  e_{n,k}=\frac{1}{2}\left(T(n,k)+(-1)^{n-k}\binom{n-1}{k-1}\right),\ \ \
  o_{n,k}=\frac{1}{2}\left(T(n,k)-(-1)^{n-k}\binom{n-1}{k-1}\right).
\end{equation*}
Let, for example, $k=2$. Then 
\begin{equation*}
  e_{n,2}=\frac{1}{2}\left(T(n,2)+(-1)^n(n-1)\right),\ \ \ 
  o_{n,2}=\frac{1}{2}\left(T(n,2)-(-1)^n(n-1)\right)
\end{equation*}
and we get the following table:
$$
\begin{array}{c|c|c|c|c|c|c|c|c|c|c|c}
 n&1&2&3&4&5&6&7&8&9&10&\dots\\
 \hline
 T(n,2)&0&1&4&11&26&57&120&247&502&1013&\dots\\
 \hline
 c(n,2)&0&1&-2&3&-4&5&-6&7&-8&9&\dots\\
 \hline
 e_{n,2}&0&1&1&7&11&31&57&127&247&511&\dots\\
 \hline
 o_{n,2}&0&0&3&4&15&26&63&120&255&502&\dots
\end{array}
$$
Let us write out all permutations of order $4$ with $2$ weak excedances:
$\underline{\mathbf{14}23}$, $\underline{\mathbf{2}1\mathbf{4}3}$, $\mathbf{24}13$,
$\underline{\mathbf{3}12\mathbf{4}}$, $\mathbf{3}1\mathbf{4}2$,
$\underline{\mathbf{34}12}$, $\mathbf{34}21$,
$\underline{\mathbf{4}1\mathbf{3}2}$, $\underline{\mathbf{42}13}$,
$\mathbf{43}12$, $\underline{\mathbf{43}21}$ ---
total $7$ even and $4$ odd permutations (weak excedances are in bold,
even permutations are underlined). 

\vspace{2mm}
{\bf Acknowledgment.} The author is grateful to  V.S. Shevelev for useful  pointers to the literature  and comments. 

\renewcommand{\refname}{Referenses}

\end{document}